\documentclass{amsart}
\usepackage{amsmath}
\usepackage{amssymb}
\usepackage{amsthm}

\usepackage{upref}
\usepackage{hyperref}
\usepackage{paralist}

\newtheorem{theorem}{Theorem}[section]
\newtheorem{lemma}[theorem]{Lemma}
\newtheorem{claim}[theorem]{Claim}
\newtheorem{corollary}[theorem]{Corollary}
\newtheorem{fact}[theorem]{Fact}
\newtheorem*{theorema}{Theorem A}

\theoremstyle{definition}
\newtheorem{definition}[theorem]{Definition}
\newtheorem{example}[theorem]{Example}

\newcommand{\RR}{\mathbb R}
\newcommand{\llex}{<_{\text{lex}}}
\newcommand{\defeq}{=}
\newcommand{\restrict}{\upharpoonright}
\newcommand{\scht}{:}

\title[Definable linear orders in o-minimal structures]{Definable linear orders
definably embed into lexicographic orders in o-minimal structures}

\author{Janak~Ramakrishnan}

\subjclass[2010]{Primary 03C64; Secondary 06A05}
\keywords{o-minimality, linear order, lexicographic order}

\begin{document}

\begin{abstract} We classify definable linear orders in o-minimal structures
expanding groups.  For example, let $(P,\prec)$ be a linear order definable in
the real field.  Then $(P,\prec)$ embeds definably in
$(\RR^{n+1},<_{\text{lex}})$, where $<_{\text{lex}}$ is the lexicographic order
and $n$ is the o-minimal dimension of $P$.  This improves a result of Onshuus
and Steinhorn in the o-minimal group context.
\end{abstract}

\maketitle

\section{Introduction}
Linear orders are defined by a simple relation, but analyzing them can be quite
difficult.  In this paper, we consider linear orders definable in o-minimal
groups and give a complete characterization.  The study of objects definable in
o-minimal structures is an active one \cite{HOP7,HO8}, since o-minimal
structures are ``tame'' and yet can be expressive enough to define objects of
interest to a wide variety of other mathematical areas \cite{PS4}.

A recent result of Onshuus and Steinhorn implies that any definable linear order
in an o-minimal structure $M$ with elimination of imaginaries is a finite union
of definable sets, each of which definably embeds in $M^n$ for some $n$, ordered
lexicographically \cite[Cor.~5.1]{OnsSte09}.  However, this result does not say
how elements are compared across sets in the union, and so the ordering is not
fully captured by this presentation.

We present an independently-discovered characterization of such definable linear
orders that completely describes the ordering when the o-minimal structure also
defines an order-reversing injection.  Say that an ordered structure $M$ with
elimination of imaginaries and such an injection is a \emph{near-group}.  The
simplest example of an o-minimal near-group is an o-minimal group\footnote{That
is, when $M$ is o-minimal and an ordered group with the o-minimal order.  Note
that $M$ may have additional structure in the form of functions, relations, etc.
(See \cite{vdD98book} for an overview of o-minimality, and \cite{Hodges93} for
basic model theory.)} with a definable positive element.

\begin{theorema}
Let $M$ be an o-minimal near-group and let $(P,\prec)$ be an $M$-definable
linear order with $n\defeq\dim(P)$.  Then there exists an $M$-definable
embedding $g$ of $(P,\prec)$ into $(M^{2n+1},\llex)$, where $\llex$ is the
lexicographic order.  Moreover, $g$ is uniformly definable over the parameters
defining $P$ and $g(P)\subseteq M^{2n+1}$ has finite projection to each odd
coordinate.
\end{theorema}

Our characterization improves that of \cite{OnsSte09} for o-minimal near-groups
since the full order is embedded in a single lexicographic order.  This means
that the study of definable linear orders in o-minimal near-groups is just the
study of definable subsets of lexicographic orders.

Besides the result in \cite{OnsSte09}, Theorem A also resembles work done in the
general context of embedding ordered sets into lexicographic products of the
reals \cite{Fleischer60,CanInd99}.  Seen from that light, Theorem A is a
definable version of results in these papers, although the results in the
general case are only partial \cite{Fleischer61}.

\subsubsection*{Outside Applications}
The study of linear orders has also been undertaken in economics.  In
\cite{BCHIM02a}, efforts were made to classify linear orders that are not
order-embeddable in the reals.  For the class of such linear orders
interpretable in o-minimal near-groups, Theorem A gives a complete
classification.

Economists have also modeled certain preference relations (linear orders) as
lexicographic orders but involving ``tradeoffs,'' in which the relative
importance of certain variables depends on their amounts \cite{Luce78,
SpaHan95}.  Theorem A shows that when such a relation is definable in an
o-minimal near-group, as it often is, the relation reduces to an associated
lexicographic order.

We note here that the uniformity in Theorem A follows from a routine
model-theoretic compactness argument.  Also, it suffices to prove Theorem A for
$\emptyset$-definable $P$, since $M$ remains an o-minimal near-group after
naming constants.

The bound of $2n+1$ is sharp by the following:

\begin{example}
Let $M\defeq(\RR,<,+,0)$ and let $n>0$.  Let $P\defeq\{\langle
x_1,\ldots,x_{2n+1}\rangle\in M^{2n+1}\scht x_i\in\{0,1\}\text{ for $i$ odd}\}$.
Let $\prec$ be the lexicographic order on $P$.
\end{example}

There is no embedding of $(P,\prec)$ into a lexicographic order of lower
dimension, due to the lack of definable injections between $M$ and proper
subsets of $M$.  However, given appropriate maps, we have:

\begin{corollary}\label{fieldbound}
If $M$ is an o-minimal field, then in Theorem A the codomain of $g$ can be taken
to be $M^{n+1}$, with $g(P)$ having finite projection to the last coordinate.
\end{corollary}

We will give the proof of Corollary \ref{fieldbound} after that of Theorem A.

The proof of Theorem A goes by induction.  The general case requires techniques
to reduce the dimension, whereas the $1$-dimensional case is more of a proof by
taxonomy.

I would like to thank C. Steinhorn for an informative discussion when I first
learned of his results with A. Onshuus, and a later discussion that helped
clarify the direction of this paper, as well as F. Wagner for an encouraging
talk on generalizing the result in the absence of a field.

\section{Notation and preliminaries}

We write ``definable'' to mean ``$\emptyset$-definable.''  Throughout, $M$ is an
o-minimal structure and $P$ a linear order definable in $M$ with $n=\dim(P)$.

\begin{definition} For $m\ge i\ge 1$, let $\pi^m_i:M^m\to M$ be projection onto
the $i$-th coordinate.  Let $\pi^m_{\le i}:M^m\to M^i$ be the map sending $x$ to
the $i$-tuple $\langle\pi^m_1(x),\ldots,\pi^m_i(x)\rangle$, and similarly for
$\pi^m_{<i}$ and $\pi^m_{>i}$.  We just write $\pi_i$, $\pi_{\le i}$, etc, since
$m$ is always clear from context.  For $x\in M^m$, let $x_i\defeq\pi_i(x)$, and
similarly for $x_{\le i}$, $x_{<i}$, and $x_{>i}$.
\end{definition}

\begin{definition}
A function $g$ is a \emph{flex-embedding of $P$} if $g$ is an embedding of
$(P,\prec)$ into $(M^{2n+1},\llex)$, with $\pi_i(g(P))$ a finite set for $i$
odd.
\end{definition}

\begin{definition}
Let $\mathcal C$ be a cell decomposition of $M^m$, and let $B\subseteq M^m$.
Define $\mathcal C\cap B\defeq\{C\in\mathcal C\scht C\subseteq B\}$.  Say that
$\mathcal C$ is \emph{compatible} with $B$ if, for every $C\in\mathcal C$,
either $C\cap B=\emptyset$ or $C\subseteq B$.  Say that $\mathcal C$ has
\emph{good projection} if, for any $i<m$ and $C,C'\in\mathcal C$, either
$\pi_{\le i}(C)\cap\pi_{\le i}(C')=\emptyset$ or $\pi_{\le i}(C)=\pi_{\le
i}(C')$.
\end{definition}

The following is a straightforward application of cell decomposition.

\begin{fact}\label{goodproj}
Let $\mathcal C$ be a definable cell decomposition of $M^m$.  There is a
definable cell decomposition $\mathcal D$ that refines $\mathcal C$ and has good
projection.
\end{fact}

\section{$1$-dimensional definable linear orders}

For a $1$-dimensional cell $C\subseteq M^m$, the order $<$ induces an order on
$C$ in a natural way via the $p_C$ function of \cite[Ch.~3(2.7)]{vdD98book}: for
$x,y\in C$, we have $x<y$ if and only if $p_C(x)<p_C(y)$.

A version of the following lemma is folklore, due to C. Steinhorn, with variants
stated in \cite{HO8} and \cite{OnsSte09}.  We need a slightly different
statement, and so we prove it here for completeness.

\begin{lemma} \label{1dfinite}
Let $M$ be an o-minimal near-group and let $(P,\prec)$ be a definable linear
order with $\dim(P)=1$.  Then $P$ is definably isomorphic to a finite union of
disjoint cells on each of which the induced $<$ and the induced $\prec$ agree.
\end{lemma}
\begin{proof}
Let $\mathcal C$ be a cell decomposition of $P$.  Fix $C\in\mathcal C$ and let
$I\defeq p_C(C)$.  The order $\prec$ induces a linear order on $I$.  For $x\in
I$, let $G(x)\defeq\{y\in I\scht y\succ x\}$.  Let $\mathcal J$ be a cell
decomposition of $I$ such that for each $J\in\mathcal J$, for all $x\in J$ the
set $G(x)$ has the same number of infinite connected components, and the
functions defining the boundaries of these components are monotonic and
continuous.  Fix an interval $J\in\mathcal J$, and let $L=\{f_1,\ldots,f_m\}$
and $U=\{g_1,\ldots,g_m\}$ be these respectively lower- and
upper-boundary-defining functions for $x\in J$.  By uniform finiteness of
families for o-minimal structures and basic properties of linear orders, there
are only finitely many $x\in I$ with $|G(x)|$ finite, so $m>0$.

Some function in $L\cup U$ must be nonconstant on $J$, since else $G(x)$ and
$G(y)$ differ by a bounded finite number of points for $x,y\in J$, which easily
violates $\prec$ being a linear order on $J$.  We show that the functions in
$L\cup U$ ``accord'' -- that if some $f\in L$ is increasing, then no $f'\in L$
is decreasing and no $g\in U$ is increasing, and similarly for the other
possibilities.  Assume that we have $f_i\in L$ increasing and $f_j\in L$
decreasing for some $i,j\le m$.  If $x<y\in J$ with $y$ sufficiently close to
$x$, then $f_j(y)<f_j(x)<g_j(y)$.  So $(f_j(y),g_j(y))\setminus
G(x)\ne\emptyset$ and $(f_i(x),g_i(x))\setminus G(y)\ne\emptyset$, which
contradicts $G$ defining a decreasing family.  The arguments for $g_i,g_j\in U$
and for $f_i\in L$, $g_j\in U$ are similar.

This ``accord'' easily implies that $\prec$ is either increasing or decreasing
on $J$.  Let $\theta$ be a definable order-reversing injection on $M$.  Fix
distinct definable $a_J\in M$ for $J\in\mathcal J$, and let $J'$ be
$\{a_J\}\times J$ if $\prec$ is increasing on $J$, and $\{a_J\}\times\theta(J)$
otherwise.  Then $\{J'\scht J\in\mathcal J\}$ is a disjoint collection of cells
on each of which the induced $\prec$ and induced $<$ agree.  Repeating this
procedure for each $C\in\mathcal C$, we are done.
\end{proof}

We can now prove Theorem A for $1$-dimensional structures.\footnote{After
proving the result, we were informed by C. Steinhorn that he already had a
version of it earlier.}  We use Lemma \ref{1dfinite} to break up a definable
linear order into cells, on each of which the order and the structure's order
agree.  We must then analyze how these pieces fit together in the definable
order.

\begin{theorem}
Let $M$ be an o-minimal near-group, and let $(P,\prec)$ be a definable linear
order with $\dim(P)=1$.  Then there exists definable $g$, a flex-embedding of
$P$.
\end{theorem}
\begin{proof}
Say that a lexicographically ordered subset of $M^3$ with finite projections to
the first and third coordinates is ``nice.''

By Lemma \ref{1dfinite}, we can suppose that $P$ has a cell decomposition,
$\mathcal D$, such that $\prec$ is increasing on each $D\in\mathcal D$ with
respect to the induced $<$.  Let $k\defeq|\mathcal D|$.  We show the theorem by
induction on $k$.  The case $k=1$ is trivial, since the unique cell
$D\in\mathcal D$ maps via the $p_D$ function into $M$ and then to $\{0\}\times
M\times\{0\}$.  We prove the case $k$, given it for case $k-1$.  Fix
$D\in\mathcal D$ and replace $P$ by $P'\cup I$, where
\begin{inparaenum}[(i)]
\item $P'$ is a nice subset of $M^3$ that is the image of the embedding of
$\mathcal D\setminus\{D\}$ into $M^3$ given by induction,
\item $I\defeq p_D(D)$ is a point or an interval in $M$, and
\item $\prec$ and $<$ agree on $I$.
\end{inparaenum}  Our concern is how $I$ and $P'$ interact.

\begin{claim}\label{samecut}
There is a cell decomposition of $I$ such that for each cell $C$ in the
decomposition, one of the following holds:
\renewcommand{\theenumi}{P\Roman{enumi}} \renewcommand{\labelenumi}{(\theenumi)}
\begin{enumerate}
\item\label{condpred} For each $x\in C$ there is $y\in P'$ with $y$ the
immediate $\prec$-successor of $x$ (that is, $x\prec y$ and
$(x,y)_\prec=\emptyset$);
\item\label{condsucc} For each $x\in C$ there is $y\in P'$ with $y$ the
immediate $\prec$-predecessor of $x$;
\item\label{condcut} Every element of $C$ lies in the same $\prec$-cut in $P'$.
\end{enumerate}
\renewcommand{\theenumi}{\arabic{enumi}}
\renewcommand{\labelenumi}{(\theenumi)}
\end{claim}

\begin{proof}
We consider a cell decomposition of $I$ compatible with the subsets defined by
the following conditions on a point $x\in I$, and show that this cell
decomposition will satisfy the claim after finitely many subdivisions.

\subsubsection*{Conditions}
\renewcommand{\theenumi}{C\arabic{enumi}}
\renewcommand{\labelenumi}{(\theenumi)}
\begin{enumerate}
\item\label{pred} there exists $y\in P'$ the immediate $\prec$-successor of $x$;
\item\label{succ} there exists $y\in P'$ the immediate $\prec$-predecessor of
$x$;
\item\label{pre} there exists $y\in P'$ with $x\prec y$ and $(x,y)_\prec\cap
P'=\emptyset$ but \eqref{pred} fails;
\item\label{suc} there exists $y\in P'$ with $y\prec x$ and $(y,x)_\prec\cap
P'=\emptyset$ but \eqref{succ} fails;
\item\label{inf} $x\succ P'$ or $x\prec P'$.
\end{enumerate}
\renewcommand{\theenumi}{\arabic{enumi}}
\renewcommand{\labelenumi}{(\theenumi)}

For any $B\subseteq I$, say that $B$ \emph{satisfies} one of the above
conditions if that condition holds for all $x\in B$.  Let $\mathcal C$ be a cell
decomposition of $I$ compatible with the sets defined by these conditions.

\begin{claim}\label{preconst}
Let $B\subseteq I$ be an interval satisfying \eqref{pre}.  Then $B$ realizes
finitely many $\prec$-cuts in $P'$.  Likewise if $B$ satisfies \eqref{suc}.
\end{claim}
\begin{proof}
For $x\in B$, let $f(x)$ denote the (necessarily unique) element of $P'$ with
$(x,f(x))\cap P'=\emptyset$.  Assume the claim fails.  Since $f(x)$ determines
the $\prec$-cut of $x$, the set $f(B)$ is infinite.  By a routine dimension
argument on fibers, there are infinitely many $x\in B$ such that the set
$f^{-1}(f(x))$ is finite.  Choose $a$ with $f^{-1}(f(a))$ finite, and in
addition with $a=\max(f^{-1}(f(a)))$.  Let $z>a$ in $B$, so $f(z)\ne f(a)$.
Since $(a,f(a))_\prec\cap P'=\emptyset$, we have $f(z)\notin (a,f(a))_\prec$.
Also, since $(z,f(z))_\prec\cap P'=\emptyset$, we have
$f(a)\notin(z,f(z))_\prec$.  Thus $z\succ f(a)$.  For any $y\in I$ with $y\succ
a$, there is $z\in B$ with $z\in(a,y)_\prec$, since $B$ is an interval.  Thus,
for all $y\in I$, we have $y\notin(a,f(a))$, so $a$ satisfies \eqref{pred},
contradiction.  The argument for \eqref{suc} is similar.
\end{proof}

\begin{claim}\label{nocond}
Let $B\subseteq I$ be a definable set such that no condition holds on $x$, for
all $x\in B$.  Then $B$ realizes finitely many $\prec$-cuts in $P'$.
\end{claim}
\begin{proof}
Let $h_1(x)\defeq\max\{\pi_1(y)\scht y\in P',y\prec x\}$.  Since $P'$ is nice,
$\pi_1(P')$ is a finite set, and so $h_1(x)$ takes only finitely many possible
values for $x\in B$.  Partitioning $B$, we suppose that $h_1(x)$ is constant on
$B$, given by $c$.  Furthermore, we suppose that for every $x\in B$, there is
$y\in P'$ with $\pi_1(y)=c$ and $y\succ x$, since the set of $x\in B$ for which
such a $y$ does not exist is definable, and all such $x$ lie in the same
$\prec$-cut of $P'$.

Let $h_2(x)\defeq\sup\{\pi_2(y)\scht y\in P',y\prec x,\pi_1(y)=c\}$.  By the
``furthermore'' supposition, $h_2(x)\in M$ for $x\in B$.  Assume that for some
$a\in B$, there exists $b\in M$ with $\langle c,h_2(a),b\rangle\in P'$.  By
niceness of $P'$, there are only finitely many $y\in P'$ with $\pi_{\le
2}(y)=\langle c,h_2(a)\rangle$.  But then $a$ must satisfy one of
\eqref{pred}-\eqref{suc} with some such $y$, contradiction.  Thus $P'$ contains
no elements with first two coordinates $\langle c,h_2(x)\rangle$ for any $x\in
B$, so in $P'$, the $\prec$-cut of $x\in B$ and the $\llex$-cut of $\langle
c,h_2(x),0\rangle$ are the same.

We can then show that the elements of $B$ realize finitely many $\prec$-cuts in
$P'$ by showing that the set $h_2(B)$ is finite.  If $h_2(B)$ were infinite then
it would contain an interval, but this is impossible, since $\langle
c,h_2(x)\rangle\notin\pi_{\le 2}(P')$ for any $x\in B$, and $h_2$ is defined as
a $\sup$ of elements in $\pi_2(P')$.
\end{proof}

This proves Claim \ref{samecut}, since if $C\in\mathcal C$ satisfies
\eqref{inf}, we can partition $\mathcal C$ so that every element lies in the
same $\prec$-cut of $P'$, and, due to Claims \ref{preconst} and \ref{nocond}, we
can partition each $C\in\mathcal C$ satisfying \eqref{pre}, \eqref{suc}, or
satisfying no conditions, so that all elements lie in the same $\prec$-cut.
\end{proof}

Fix a cell decomposition of $I$ satisfying Claim \ref{samecut},
$I_1<\cdots<I_m$.  Note that $\{I_2,\ldots,I_m\}$ is a cell decomposition of
$I\setminus I_1$ satisfying Claim \ref{samecut} with respect to $P'\cup I_1$,
since properties \eqref{condpred} and \eqref{condsucc} are trivially preserved,
and $I_1\prec I\setminus I_1$ implies that property \eqref{condcut} is too.

We will give a definable embedding $g$ of $P'\cup I_1$ into $M^3$ such that the
image is still nice.  The decomposition $I_2,\ldots,I_m$ will satisfy Claim
\ref{samecut} with respect to $g(P'\cup I_1)$ by the above argument, so we will
be done by induction on $m$.

If $I_1$ satisfies property \eqref{condpred} of Claim \ref{samecut}, let
$f:I_1\to P'$ be the definable injection with $f(x)$ the unique $y\in P'$ such
that $y\succ x$ and $(x,y)_\prec$ is empty.  Let $f_3(x)\defeq\pi_3(f(x))$.  By
niceness of $P'$, for any $x\in I_1$ the set $R_x\defeq\{z\prec f(x)\scht
\pi_{\le 2}(f(x))=\pi_{\le 2}(z)\}$ is finite.  Thus, by elimination of
imaginaries there is some definable function $h$ with $h(x)<f_3(x)$ and
$\langle\pi_{\le 2}(f(x)),h(x)\rangle\succ R_x$.  Let $g(x)\defeq\langle
\pi_{\le 2}(f(x)),h(x)\rangle$ for $x\in I$ and extend $g$ on $P'$ by the
identity.  The function $g$ is a definable embedding of the ordered set $P'\cup
I_1$ into $M^3$ ordered lexicographically.  For $x\in I_1$, we have $\pi_{\le
2}(f(x))=\pi_{\le 2}(f(y))$ for only finitely many $y\in I_1$, which implies
that $g(P'\cup I_1)$ is nice.

We proceed analogously if $I_1$ satisfies property \eqref{condpred} with respect
to $P'$.

Now suppose that $I_1$ satisfies property \eqref{condcut} with respect to $P'$.
First, suppose that this cut is also satisfied by some $\langle a,0,0\rangle$
with $a\notin\pi_1(P')$.  Then map $I_1$ to $\langle a,I_1,0\rangle$, and fix
$P'$.  It is easy to verify that this map has the desired properties.
Otherwise, there are $b,b'\in P'$ with $\pi_1(b)=\pi_1(b')$ and $b\prec I_1\prec
b'$.  Let $a\defeq\pi_1(b)$.  Let $B\defeq\{x\in P'\scht\pi_1(x)=a\land x\succ
I_1\}$.  Let $c$ be the least element of $\pi_1(P')$ greater than $a$, or
$\infty$ if $a=\max(\pi_1(P'))$.  Choose definable elements $d,e\in M$ with
$a<d<e<c$.  Let $g:P'\cup I_1\to M^3$ be the identity on $P'\setminus B$, and
let $g$ send $x\in B$ to $\langle e,\pi_{\ge 2}(x)\rangle$ and send $I_1$ to
$\langle d,I_1,0\rangle$.  The map $g$ is a definable embedding of $P'\cup I_1$
into $M^3$.  It is easy to see that $g(P'\cup I_1)$ is still nice.
\end{proof}

\section{$n$-dimensional definable linear orders}

\begin{proof}[Proof of Theorem A]
Let $H\defeq\{x\in P\scht\forall y\prec x(\dim((y,x)_\prec)=n)\}$.  Points in
$H$ have ``intrinsically full dimension'' below -- any $\prec$-interval
approaching one from below has dimension $n$.  Note that $H$ is definable (see
\cite[Ch.~4(1.5)]{vdD98book}).

\begin{lemma}\label{pdim=n}
If $\dim(H)=n$, then $n\le 1$.
\end{lemma}
\begin{proof}
For each $x\in H$, let $B_x\defeq\{z\in P\scht (x,z)_\prec\text{
infinite},(x,z)_\prec\cap H=\emptyset\}$.  For distinct $x,y\in H$, the sets
$B_x$ and $B_y$ are disjoint, and if $B_x$ is nonempty, it has positive
dimension.  Thus, the set $B\defeq\{x\in H\scht B_x\ne\emptyset\}$ must have
dimension less than $n$.

Let $\Gamma\subseteq H\setminus B$ be a definable connected $1$-dimensional set.
Applying Lemma \ref{1dfinite} and restricting, we may suppose that $\prec$
agrees with the induced $<$ on $\Gamma$.  Let $T:P\to\Gamma$ be the definable
partial function $T(x)\defeq\inf_{\prec\restrict\Gamma}\{y\in\Gamma\scht
y\succeq x\}$ when this $\inf$ exists.  Note that $T(x)$ is defined if there is
$z\in\Gamma$ with $z\prec x$.  The sets $T^{-1}(y)$ and $T^{-1}(z)$ are
$\prec$-convex, and disjoint for distinct $y,z\in\Gamma$.  By cell
decomposition, there is a definable infinite connected set
$\Gamma'\subseteq\Gamma$ on which $\dim(T^{-1}(y))$ is constant.  Let
$p\defeq\dim(T^{-1}(y))$ for $y\in\Gamma'$.

Fix $b,c\in\Gamma'$ with $b\prec c$.  If $x\in(b,c)_\prec$ then $b\preceq
T(x)\preceq c$.  Therefore
$(b,c)_\prec\subseteq\bigcup_{y\in\Gamma'\cap[b,c]_\prec} T^{-1}(y)$, and so
$\dim((b,c)_\prec)\le\dim\left(\bigcup_{y\in\Gamma'\cap[b,c]_\prec}
T^{-1}(y)\right)=1+p\le n$.

Assume for a contradiction that $p>0$.  Fix $a\in\Gamma'$.  First, assume that
there exists $d\in T^{-1}(a)$ with $d\prec a$.  Then $(d,a)_\prec\subseteq
T^{-1}(a)$, so $\dim((d,a)_\prec)\le p<n$, but since $a\in H$, we have
$\dim((y,a)_\prec)=n$ for all $y\prec a$, contradiction.  Thus $T^{-1}(a)\succeq
a$.  Fix $d\in T^{-1}(a)$ with $(a,d)_\prec$ infinite -- possible since
$\dim(T^{-1}(a))>0$.  Since $a\in H\setminus B$, there is $d'\in
H\cap(a,d)_\prec$, so $\dim((a,d')_\prec)=n\le p<n$, contradiction.

Thus $p=0$ and $\dim((b,c)_\prec)\le 1$, so $n\le 1$.
\end{proof}

If Lemma \ref{pdim=n} holds, then we are done by the $1$-dimensional case, so we
suppose from now on that $\dim(H)<n$.

Let $E$ be the equivalence relation on $P$ defined as $xEy$ if and only if
$\dim((x,y)_\prec\cup (y,x)_\prec)<n$.  Note that the $E$-classes of $P$ are
$\prec$-convex.

\begin{lemma}\label{infclasses}
No $E$-class has dimension $n$.
\end{lemma}
\begin{proof}
Assume not, so there is an $E$-class $B$ with $\dim(B)=n$.  We replace $P$ by
$B$.  Then for any $x,y\in P$, $\dim((x,y)_\prec)<n$, but $\dim(P)=n$.  Consider
the partial $M$-types $p_1(x)$, which says that $x\in P$ and $x\prec a$ for each
$a\in P(M)$; and $p_2(x)$, which says that $x\in P$ and $x\succ a$ for each
$a\in P(M)$, and let $b_i\models p_i$ for $i=1,2$ with $b_1,b_2\in M'$, an
elementary extension of $M$.  Then $\dim((b_1,b_2)_\prec)<n$, since this
first-order property of $P$ is preserved in $M'$, but $P(M)$ is $n$-dimensional
and contained in $(b_1,b_2)_\prec$, contradiction.
\end{proof}

Lemma \ref{infclasses} implies that $E$ has infinitely many equivalence classes.
The proof now proceeds through quotienting by $E$.  We first show that
$m\defeq\dim(P/E)<n$.  If not, then there is $B\subseteq P/E$ with $\dim(B)=n$
such that each $E$-class represented in $B$ is finite.  Each $E$-class
represented in $B$ has a $\prec$-least element.  Let $D$ be the set of these
$\prec$-least elements, so $\dim(D)=\dim(B)=n$.  For any $x\in D$ and any
$y\prec x$, we have $\dim((y,x)_\prec)=n$, so if $x\in D$ is not the
$\prec$-least element of $P$, then $x\in H$.  Thus $\dim(D)=n$ implies
$\dim(H)=n$, contradiction.

The order on $P$ induces a linear order on $P'\defeq P/E$.  By induction, there
exists $g$ a definable embedding of $(P',\prec)$ into $(M^{2m+1},\llex)$, with
$g(P')$ having finite projection to each odd coordinate.  Then $P$ is definably
isomorphic as an order to $\{\langle x,y\rangle\scht x\in g(P'),y\in
[g^{-1}(x)]_E\}$, ordered lexicographically by the orders $\llex$ on $g(P')$ and
$\prec$ on $[g^{-1}(x)]_E$ for $x\in P'$.  We replace $P$ by this ordered set
and $P'$ by $g(P')$.  Let $Q_x\defeq[x]_E$, ordered by $\prec$, for $x\in P'$.
The remainder of the proof is just to bound the dimension of the embedding of
$P$, since it is easy to embed the $Q_x$'s uniformly in some lexicographic
order.

\subsubsection*{Compressing $P'$}

Let $\mathcal C$ be a cell decomposition of $P'$ with good projection
such that, for $C\in\mathcal C$, if $x,y\in C$ then $\dim(Q_x)=\dim(Q_y)$.

We must ``compress'' each $C\in\mathcal C$ while preserving the lexicographic
order.  For odd $j<2m$, $i\in(j,2m+1]$, and $C\in\mathcal C$, let
\begin{equation*}
V_i(j,C)\defeq\{D\in\pi_{\le i}(\mathcal C)\scht
\pi_{<i}(D)=\pi_{<i}(C),D\ne\pi_{\le i}(C)\}.
\end{equation*}
Let $k(j,C)$ be the greatest coordinate $k$ such that $\dim(\pi_{\le
k}(C))=\dim(\pi_{\le j}(C))$, and let $V(j,C)=\bigcup_{j<i\le k(j,C)}V_i(j,C)$.
The collection $V(j,C)$ represents cells that must be shifted before the
coordinates of $C$ can be collapsed.  If $V(j,C)=\emptyset$, then transforming
$C$ by modifying coordinates $>j$ will not affect the ordering between elements
of $C$ and the rest of $P'$.

Fix odd $j$ minimal and $C\in\mathcal C$ such that $V(j,C)\ne\emptyset$, and let
$k\defeq k(j,D)$.  For each $i\in(j,k]$, let $r(i)\defeq|V_i(j,C)|$.  Each
$V_i(j,C)\cup\{\pi_{\le i}(C)\}$ is totally ordered by the relation $D<_iD'$
given by
\begin{equation*}
\forall x\in\pi_{<i}(C)\left(\{y\scht \langle x,y\rangle\in D\}<\{y\scht \langle
x,y\rangle\in D'\}\right).
\end{equation*}
Let $s(i)\defeq\left|\left\{D\in V_i(j,C)\scht D<_i C\right\}\right|$.  Let
$\pi_j(C)=\{b\}$.  For $i=j+1,\ldots,k$, fix definable $c^i_1<\cdots<c^i_{r(i)}$
in $M$ such that the following three properties hold:
$c^i_{s(i)}<b<c^i_{s(i)+1}$; for $i>j+1$ we have
$c^{i-1}_{s(i-1)}<c^i_1<c^i_{r(i)}<c^{i-1}_{s(i-1)+1}$; and
$(c^{j+1}_1,c^{j+1}_{r(i)})\cap \pi_j(P')=\{b\}$.  The $c^i$'s define nested
intervals around $b$.  Then let $h:P'\to M^{2m+1}$ be defined by $h(x)\defeq x$
if $\pi_{\le i}(x)\notin \bigcup V(j,C)$ for all $i>j$, and otherwise,
$h(x)\defeq\langle x_{<j},c^i_t,x_{>j}\rangle$, where $\pi_{\le i}(x)$ belongs
to the $t$-th member of the $<_i$-ordered set $V_i(j,C)$.  The function $h$ maps
the finitely many ways at which a cell can ``branch'' from $C$, at coordinates
$j+1,\ldots,k$, into the $j$-th coordinate.

The function $h$ is an embedding of $(P',\llex)$ and the set $h(\mathcal C)$ is
a cell decomposition of $h(P')$ with good projection.  Moreover, if we let $V_h$
be defined as $V$ was but for $h(P')$, then for any $C'\in h(\mathcal C)$ and
any odd $j'<2m$, we have $|V_h(j',h(C'))|\le|V(j',C')|$, and in particular
$V_h(j,h(C))=\emptyset$.  Thus, after replacing $P'$ by $h(P')$ and repeating
this finitely many times, we may suppose that $V(j,C)=\emptyset$ for all odd $j$
and all $C\in\mathcal C$.

Fix $C\in\mathcal C$ and let $k(C)$ be the first odd coordinate $k$ such that
$\dim(\pi_{\le k+1}(C))=\dim(\pi_{\le k}(C))$, or $2m+1$ otherwise.  Let
$O\defeq\{j\in(k(C),2m]\scht\dim(\pi_{\le j}(C))>\dim(\pi_{<j}(C)\}$, and let
$O=\{j(1),\ldots,j(r)\}$.  Define $h_C:C\to M^{2m+1}$ by $h_C(x)=\langle x_{\le
k(C)},x_{j(1)},0,x_{j(2)},0,\ldots,x_{j(r)},0,\ldots,0\rangle$, and let $h$ be
the union of the $h_C$'s.  We show that $h$ is an order-preserving embedding of
$P'$.  For distinct $C_1,C_2\in\mathcal C$, let $s\defeq\min(k(C_1),k(C_2))$.
Since $V(k(C_i),C_i)=\emptyset$ for $i=1,2$, if $x\in C_1$ and $y\in C_2$ then
$\pi_{\le s}(x)\ne\pi_{\le s}(y)$, and since $\pi_{\le s}(h(x))=\pi_{\le s}(x)$,
and similarly for $y$, the map $h$ must preserve the ordering on $x$ and $y$.
Thus, we can restrict to a single $C$.  Let $x,y\in C$.  Let $i$ be the first
coordinate such that $x_i\ne y_i$.  If $i<k(C)$, then $\pi_{\le
i}(h(x))=\pi_{\le i}(x)$ and similarly for $y$ and we are done.  Thus $i\in O$,
so $i=j(t)$ for some $t\le r$.  By definition of $h$, the first coordinate at
which $h(x)$ and $h(y)$ differ is $l\defeq k(C)+2t-1$, at which
$\pi_l(h(x))=\pi_i(x)$ and $\pi_l(h(y))=\pi_i(y)$, so we have shown that $h$ is
an order-preserving embedding of $P'$.  Moreover, $h(\mathcal C)$ is a cell
decomposition of $h(P')$.

After replacing $P'$ by $h(P')$ and $\mathcal C$ by $h(\mathcal C)$, we have
$\pi_i(C)=\{0\}$ for all $i>2\dim(C)+1$ and $C\in\mathcal C$.

\subsubsection*{Compressing $Q_x$}
For $C\in\mathcal C$, let $q(C)=\dim(Q_x)$ for $x\in C$.  By Lemma
\ref{infclasses}, $q(C)<n$.  By induction for Theorem A, each $Q_x$ can be
definably flex-embedded in $M^{2q(C)+1}$, and since $Q_x$ is uniformly definable
in $x$, this embedding can be taken to be uniform as well, so we have a
definable function $g_C$ with $g_C(x,-)$ a flex-embedding of $Q_x$ for all $x\in
C$.  Letting $q\defeq\max\{q(C)\scht C\in\mathcal C\}$, we can embed each
$M^{2q(C)+1}$ in $M^{2q+1}$, extending by $0$, and so suppose that each $g_C$
embeds into $M^{2q+1}$, and let $g$ be their union.  Replace $P$ by $\{\langle
x,g(x,y)\rangle\scht x\in P',y\in Q_x\}$, so $Q_x$ is replaced by $g(x,Q_x)$.
For each odd $i\le 2q+1$ and $x\in P'$, the set $\pi_i(Q_x)$ is finite.  Thus,
by o-minimality $|\pi_i(Q_x)|$ is bounded as $x$ ranges over $P'$, and so we can
set $r\defeq\max\{|\pi_i(Q_x)|\scht \text{odd }i\le 2q+1,x\in P'\}$, and fix
definable $a_1<\cdots<a_r\in M$.  Then define $h_x:Q_x\to M^{2q+1}$ so that
$\pi_i(h_x(y))=y_i$ for $i$ even, and $\pi_i(h_x(y))=a_t$ for $i$ odd, with
$y_i$ the $t$-th element in the finite ordered set $\pi_i(Q_x)$.  Replace $P$ by
$\{\langle x,h_x(y)\rangle\scht x\in P',y\in Q_x\}$ and $Q_x$ by $h_x(Q_x)$.

\subsubsection*{Joining}
We now have $P\subseteq M^{2(m+q)+2}$, ordered lexicographically, with
$\pi_i(P)$ finite for odd $i\le 2m+1$ and even $i\ge 2m+2$.  Let $m(C)=\dim(C)$
for $C\in\mathcal C$.  For each $C\in\mathcal C$, we have $\pi_i(C)=\{0\}$ if
$2m(C)+1<i\le 2m+1$, and for $x$ with $\pi_{\le 2m+1}(x)\in C$, we have
$\pi_i(x)=\{0\}$ if $2m+1+2q(C)+1<i\le 2(m+q)+2$.

For each $C\in\mathcal C$, let $\pi_{2m(C)+1}(C)=\{b^C\}$.  Fix definable
$c^C_1<\cdots<c^C_r$ such that $b^C\in(c^C_1,c^C_r)$ and
$(c^C_1,c^C_r)\cap(c^D_1,c^D_r)=\emptyset$ for $D\ne C$.  Let $g_C$ take $x\in
P$ with $\pi_{\le 2m+1}(x)\in C$ to $\langle x_{\le
2m(C)},c_t,x_{2m+3},x_{2m+4},\ldots,x_{2m+2q(C)+2}\rangle$, where $x_{2m+2}$ is
the $t$-th element of the ordered set $\{y\scht\langle
x_{<2m+2},y\rangle\in\pi_{\le 2m+2}(P)\}$.  Note that this ordered set has at
most $r$ elements.  The codomain of $g_C$ is $M^{2m(C)+2q(C)+1}$.  Since
$m(C)+q(C)\le n$, we can take all the $g_C$'s to map to $M^{2n+1}$ through
extending by $0$.  Then the union of the $g_C$'s is the desired embedding.

\end{proof}

\begin{proof}[Proof of Corollary \ref{fieldbound}]

The bound in Corollary \ref{fieldbound} comes from taking the image of $g(P)$
under embeddings whose existence is guaranteed by the following:

\begin{claim}
Let $M$ be an o-minimal field.  Let $B\subset M^n$ be definable, with
$|\pi_k(B)|$ finite for some $k<n$.  Then $(B,\llex)$ embeds definably into
$(M^{n-1},\llex)$.
\end{claim}
\begin{proof}
Let $\pi_k(B)=\{a_1<\cdots<a_m\}$.  Let $a_0\defeq -\infty$.  For $i\le m$, let
$f_i:M\to(a_{i-1},a_i)$ be a definable order-preserving injection.  For $x\in
B$, let $h(x)\defeq\langle x_{<k},f_i(x_{k+1}),x_{>k+1}\rangle$, when $x_k=a_i$.
Then $h$ is the desired embedding.
\end{proof}\renewcommand{\qedsymbol}{}
\end{proof}
\bibliographystyle{alpha-elink}
\bibliography{../janak}

\end{document}